\documentclass{article} \usepackage{amsmath,amsfonts,amssymb}
\numberwithin{equation}{section}

\newtheorem{theorem}[equation]{Theorem}
\newtheorem{lemma}[equation]{Lemma}

\newtheorem{proposition}[equation]{Proposition}
\newtheorem{definition}[equation]{Definition}
\newtheorem{definitions}[equation]{Definitions}
\newtheorem{example}[equation]{Example}

\newtheorem{remark}{Remark}

\newcommand{\Pic}{\operatorname{Pic}}
\newcommand{\Cl}{\operatorname{Cl}}

\newcommand{\PP}{{\mathbb P}}
\newcommand{\ZZ}{{\mathbb Z}}

\renewcommand{\O}{{\mathcal O}}
\newcommand{\E}{{\mathcal E}}
\newcommand{\I}{{\mathcal I}}
\newcommand{\ci}{{C^\prime}}
\newcommand{\zi}{{Z^\prime}}
\newcommand{\cii}{{C^{\prime\prime}}}
\newcommand{\ciii}{{C^{\prime\prime\prime}}}
\newcommand{\ziii}{{Z^{\prime\prime\prime}}}

\newenvironment{pf}
{\noindent\textbf{Proof.}}
{\hfill{$\square$}\medskip}

{\hfill\textbf{ Q.E.D.}\medskip}

\begin{document}

\title{Curves of maximal genus in $ \PP^5$}
\author{{\sc Rita Ferraro} 
\thanks{
Dipartimento di Matematica,
Universit\`a di Roma Tre,
Largo San Leonardo Murialdo, 1- 00146 Roma.
E-mail:{\tt ferraro@mat.uniroma3.it}.}
}


\date{}
\maketitle

\section{Introduction}

Let $C\subset{ \PP}^r$ be a reduced, irreducible, not degenerate curve,
not contained on surfaces of degree $< s$; when $d=\deg C$ is large with
respect to $s$, the {\it arithmetic} genus $p_a(C)$ is bounded by a function
$G(d, r, s)$ which is of  type $\frac{{d}^2}{2s}+ O(d)$. 

The existence
of such a  bound for $C\subset { \PP}^3$ was announced by Halphen
in 1870 and proved by Gruson and Peskine in \cite{gp}; for curves in
${ \PP}^r$, $r\ge 4$ the bound is stated and proved 
in \cite{ccd} ( for $d>\frac{2s}{r-2}\Pi_{i=1}^{r-2}
{((r-1)!s)}^{\frac{1}{r-1-i}}$).

The existence of curves of maximal genus, i.e. whose genus attains the
bound, is known in ${ 
\PP}^3$ for $d>s^2-s$, in ${ \PP}^4$ for $d>12 s^2$ and in ${\PP}^r$,
$r\ge 5$, at least for $d>>s$.
\cite{gp} contains a precise description of those  curves which do
not
lie on surfaces of degree $<s$ and whose genus is maximal: they are 
arithmetically Cohen-Macaulay, lie on a surface $S$ of degree $s$ and
they are
directly linked to plane curves.
\cite{cc} contains the description of curves in ${\PP}^4$ of
maximal genus $G(d, 4, s)$.

The  complete description of curves $C\subset
{ 
\PP}^5$, not contained on surfaces of degree $<s$, whose genus is $G(d, 5, 
s)$ has been given by the author in her PhD dissertation [F1]. The main
result of this note is the classification Theorem \ref{teorema}, which
holds for $s\geq 9$. Due to the long list of cases,
some proofs are given only in some specific examples.
We already know  \cite{ccd} that such  curves must be arithmetically
Cohen-Macaulay and they must lie on a surface $S$ of degree $s$, whose
general hyperplane section $\Gamma$ is a "Castelnuovo curve" in ${
\PP}^4$, i.e. a curve in ${ \PP}^4$ of maximal genus. When $s\geq 9$ the
surface $S$ lies
on
a
rational normal $3$-fold $X$ of degree $3$ in $\PP^5$, which can be
singular.
Analogously to \cite{gp}
and  \cite{cc},
we
describe our curves $C$ of genus $G(d, 5, s)$ in terms of the curve
$C^\prime$
obtained by linking $C$ with $S$ and a hypersurface $F$ of minimal degree
passing through $C$ and not containing $S$.

In \cite{dc}, the author derives an upper bound for the {\it geometric}
genus of integral curves on the three dimensional
{\it nonsingular} quadric which lie 
on an integral surface of degree $2k$, as a function of $k$ and the degree
$d$ of the curve, without any assumption on the degree $d$.
The author analyzes the curves that should achieve that bound, which turns out
to be {\it not sharp}.
\medskip

We state now the main Theorem. In Propositions \ref{plane}, \ref{quadric}
and \ref{cubic}
we will give a closer
description of cases 2) 3) and 4) of the Theorem.
\begin{theorem}
\label{teorema}
Let $C\subset \PP^5$ be an integral non degenerate curve not contained 
on surfaces of degree $<s$ and let $s\geq 9$; put $d=\deg C$ and let
$p_a(C)=G(d, 5, s)$ be
its arithmetic genus. Assume  
$d>\frac{2s}{3}\Pi_{i=1}^3 {(4!s)}^{\frac{1}{4-i}}$. 

Then $C$ is
arithmetically Cohen-Macaulay and lies on a irreducible
surface $S$ of degree $s$ contained in cubic a rational normal  $3$-fold
$X\subset
\PP^5$.

Put $d-1=sm+\epsilon$, $0\le \epsilon\le s-1$ and $s-1=3w+v$, $v=0, 1, 2$;

if $\epsilon< w(4-v)$, divide $\epsilon=kw+\delta$, $0\le\delta<w$;

if $\epsilon \ge w(4-v)$, divide $\epsilon+3-v=k(w+1)+\delta$, $0\le
\delta<w+1$.

Then $C$ is contained on a hypersurface $F$ of degree $m+1$, not passing
through $S$. If $C^\prime$ is the curve linked to $C$ by $F$ and $S$, we
have:
\begin{enumerate}
\item 
when $k=3$, then $C^\prime=\emptyset$; i.e. $C$ is complete
intersection on $S$;
\item 
when $k=2$, then $C^\prime$ is a plane curve;
\item 
when $k=1$, then  $C^\prime$ lies on a
surface
of degree $2$;
\item
when $k=0$, then
$C^\prime$
lies on a surface
of degree $3$.
\end{enumerate}
Finally, for all $s, d$ with  $s\geq 4$ and
$d>\frac{2s}{3}\Pi_{i=1}^3 {(4!s)}^{\frac{1}{4-i}}$, one can find a smooth 
curve
of degree $d$, arithmetic genus $G(d, 5, s)$ which does not lie on 
surfaces of degree $< s$.

\end{theorem}

The proof is based on the analysis of the Hilbert function of a general
hyperplane section $Z$ of $C$.
The main technical problem that one does not find in the previous
cases
($r= 3, 4$) is that for describing
$C^\prime$ we have to perform a linkage by a complete intersection on the
scroll $X$, which is in general, if $X$ is singular, a non-Gorenstein
scheme. To this purpose the author has proved in \cite {f2} and \cite{f3}  
 some general results 
to which we
will refer in these note. 

The last section is devoted to examples. We produce smooth curves of
maximal
genus for all $s\geq 4$ and $d>\frac{2s}{3}\Pi_{i=1}^3
{(4!s)}^{\frac{1}{4-i}}$.
It should be observed that in \cite{cc} the authors  don't analyze the
regularity
of the produced extremal curves in $\PP^4$
and that in \cite{ccd} the produced examples of curves
of maximal genus $G(d, r, s)$ in $\PP^r$ for $d\gg s$ are in general
singular. 

With the same techniques used for the
classification in $\PP^5$ it is possible to classify curves in $\PP^r$ of
maximal
genus $G(d, r, s)$ for every $r$ and $s\geq 2r-1$. In \cite {f3} the
author has given an
example of the classification procedure for curves of maximal genus
$G(d, r, s)$ in $\PP^r$ and of the costruction of such smooth extremal
curves.

\medskip\noindent
{\bf Acknowledgments}
The paper has been written while the author was supported by a Post-Doc
scholarship of Universit\'a di Roma Tre.
The author thanks Ciro Ciliberto, for his patience and help.

\section{ Weil divisors on $X$}

We will see in the next section that curves $C\subset\PP^5$ of maximal
genus $G(d, 5,s)$ which we want to
classify lie on a cubic rational normal $3$-fold $X\subset \PP^5$ about which
 we
need to fix some notation and mention some result.
A rational normal $3$-fold $X\subset\PP^5$ is the image
of a projective bundle $\pi:\PP({\mathcal\E})\to \PP^1$ over
$\PP^1$, 
via
the morphism
$j$
defined by the tautological line bundle $\O_{\PP(\E)}(1)$,
where $\E$ is a locally free sheaf 
of rank $3$ on $\PP^1$ of one of the following three kinds:
\begin{enumerate}
\item
$\mathcal E=\O_{\PP^1}(1)\oplus\O_{\PP^1}(1)\oplus\O_{\PP^1}(1)$.
In this case  $X\cong\PP^1\times\PP^1$ is  smooth and it is ruled by
$\infty^1$ disjoint
planes; we put $X=S(1, 1, 1)$. 
\item
$\mathcal E=\O_{\PP^1}\oplus\O_{\PP^1}(1)\oplus\O_{\PP^1}(2)$.
Here $X$ is a cone
over a smooth cubic surface in $\PP^4$ with vertex a point $V$
and it is ruled by $\infty^1$ 
planes intersecting at $V$; we put $X=S(0, 1, 2)$. 
\item
$\mathcal E=\O_{\PP^1}\oplus\O_{\PP^1}\oplus\O_{\PP^1}(3)$.
$X$ is a cone
over a twisted cubic in $\PP^3$ with vertex a line $l$,
and it is ruled by
$\infty^1$ 
planes intersecting at $l$; we put $X=S(0, 0, 3)$.
\end{enumerate}

Let us denote $\PP(\E)=\tilde X$. The morphism $j:\tilde X\to X$ is a
rational resolution of singularities, called the {\it canonical
resolution} of $X$. It is well known that the Picard group $\Pic (\tilde
X)$ of $\tilde X$ is isomorphic to $\ZZ[\tilde H] \oplus \ZZ[\tilde R]$,
where
$[\tilde H]=[\O_{\tilde X}(1)]$ is the hyperplane class and 
$[\tilde R]=[\pi^*\O_{\PP^1}(1)]$ is the class of the fiber of the map
$\pi:\tilde X\to \PP^1$. The
intersection form on $\tilde X$ is determined by the rules:
\begin{equation}
\label{intform}
{\tilde H}^3=3  
\qquad {\tilde R}\cdot {\tilde H}^2=1 
\qquad {\tilde R}^2\cdot {\tilde H}=0. 
\end{equation}

The cohomology of the invertible sheaf $\O_{\tilde X}(a\tilde H +
b\tilde R)$ associated to a divisor $\sim a\tilde H+b\tilde R$ in
$\tilde X$ can be explicitly
calculated using
the Leray spectral sequence. In particular,
for $a\geq 0$ and $b\geq -1$, the dimension $h^0(\O_{\tilde X}(a\tilde H+ 
b\tilde R))$ 
does not depend on the type of the scroll and it is given
by the formula (\cite{f2} 3.5):
\begin{equation}
\label{accazero}
h^0(\O_{\tilde X}(a\tilde H + b\tilde R ))
=3{{a+2}\choose{3}}+(b+1){{a+2}\choose{2}}.
\end{equation} 
Let $H$ and $R$ be the {\it strict images} in $X$ of $\tilde H$ and
$\tilde R$
respectively, i.e. the scheme-theoretic closures 
$\overline{j({\tilde H}_{|j^{-1}X_S})}$
and $\overline{j({\tilde R}_{|j^{-1}X_S})}$,
 where $X_S$ denotes the smooth part of $X$.
Let us consider on $X$ the direct image
of $\O_{\tilde X}(a\tilde H +b\tilde R)$
in $X$ for every $a, b\in \ZZ$:
\[
\O_X(aH+ bR):=j_*\O_{\tilde X}(a\tilde H + b\tilde R).
\]
If the scroll $X$ is smooth, then the sheaves $\O_X(aH+bR)$ are the invertible
sheaves associated to the Cartier divisors $\sim aH+bR$, while when $X$
is singular this is no longer true. In this case 
we have the following Proposition which has been proved in \cite{f2} (Lemma
2.14, Cor. 3.10 and Th. 3.17). The reader may refer to \cite{h1} for a survey
on {\it divisorial sheaves} associated to generalized divisors, in particular
to Weil divisors.

\begin{proposition}
\label{weil}
Let $X\subset \PP^5$  be a rational normal $3$-fold and let $j:\tilde
X\to X$ be its canonical resolution. Let $\Cl(X)$ be the group of Weil
divisors on $X$ modulo linear equivalence. Then
\begin{itemize}
\item
If $X$ is smooth or $X=S(0, 1, 2)$ we have 
 $\Cl (X)\cong\ZZ[H]\oplus \ZZ [R]$. The divisorial sheaf
associated to a divisor $\sim aH+bR$ on $X$ is $\O_X(aH+ bR)$ for every
$a, b\in \ZZ$.
\item
If $X=S(0, 0, 3)$ we have that $H\sim 3R$ and $\Cl(X)\cong \ZZ[R]$. The
sheaves $\O_X(aH +bR)$ with $a+3b=d$ and $b<3$ are all isomorphic to the divisorial
sheaf associated to a divisor $\sim dR$.   
\end{itemize}
\end{proposition}
\begin{remark}
\label{dimD}
Since $R^ij_*\O_{\tilde X}(a\tilde H+ b\tilde R)=0$ for $i>0$ and for
all $a\in 
\ZZ$ and $b\geq -1$, then formula (\ref{accazero}) holds 
for $h^0(\O_X(aH+bR))$ too.
In particular  we can use formula (\ref{accazero}) to compute $\dim|D|$ for an effective
divisor $D\sim aH+bR$ in $X=S(1, 1, 1)$ or $X=S(0, 1, 2)$ with $b\geq -1$.
When $X=S(0, 0, 3)$ by Prop. \ref {weil} we can write the divisorial sheaf associated to 
$D\sim dR$ in the form $\O_X(aH+bR)$ with $0\leq b<3$; therefore we can use formula 
(\ref{accazero}) to compute $\dim|D|$ for every effective divisor $D$ in $X=S(0, 0, 3)$.
\end{remark}

When $X$ is smooth or $X=S(0, 1, 2)$ the intersection form (\ref{intform}) 
on $\tilde X$
determines via the isomorphism $\Pic(\tilde X)\cong \Cl(X)$ of Prop.
\ref{weil}  the intersection form on $X$. In particular we use
the intersection number $D\cdot D^\prime\cdot H$ 
of two effective divisor $D\sim aH+bR$ and $D^\prime \sim a^\prime 
H+b^\prime R$ with no common components to
compute the degree of their scheme-theoretic intersection $D\cap
D^\prime$.

When $X=S(0, 0, 3)$ the computation of $\deg (D\cap D^\prime)$
is more complicate and it is explained in detail in  \cite{f2};
 here we want just to state the results that we need. So
let us first introduce the {\it integral total tranform} of an effective 
Weil divisor
$D\sim dR$:
\begin{definition}
\label{tot}
Let $X=S(0, 0, 3)$. Let $D\subset X$ be an effective Weil divisor. Then
the
integral total transform $D^*$ of $D$ in $\tilde X$ is:
\[
D^*:= \tilde D +\lceil q \rceil E
\]
where $\tilde D\sim a \tilde H+b\tilde R$ is the proper transform of $D$
in $\tilde X$, 
$E\sim \tilde H-3\tilde R$ is the exceptional divisor in $\tilde X$ and 
$\lceil q\rceil$ is the smallest integer $\geq q:=\frac{b}{3}$.
\end{definition}
Then let us introduce, for every effective Weil divisor $D$ in $X$,
the rational number $\epsilon:= \lceil q \rceil -q$.
 We can compute the degree of the intersection scheme $D\cap D^\prime$
of two effective divisors $D\sim dR$ and $D^\prime\sim d^\prime R$ 
with no common components
using the following formula which has been proved in \cite{f2} Prop. 4.11:
\begin{equation}
\label{degree}
\deg(D\cap D^\prime)=
\begin{cases}
D^*\cdot {D^\prime}^*\cdot{\tilde H}\qquad\text{if
$[\epsilon+\epsilon^\prime]=0$}\\
D^*\cdot {D^\prime}^*\cdot{\tilde H}+3(\epsilon+\epsilon^\prime-1)+1\qquad
\text{if $[\epsilon+\epsilon^\prime]=1$}.\\
\end{cases}
\end{equation}
By abusing notation we will write the degree $\deg (D\cap D^\prime)$
by the intersection number $D\cdot D^\prime\cdot H$.
Moreover we compute the intersection multiplicity 
$m(D, D^\prime;l)$ of $D$ and $D^\prime$
along the singular line $l$ of $X$ as follows:
\begin{equation}
\label{intmult}
m(D, D^\prime;l)= D\cdot D^\prime \cdot H-{\tilde D}\cdot{{\tilde D}^\prime}\cdot \tilde H.
\end{equation}

\begin{remark}
\label{description1}
To prepare the proof of the main Theorem and of Propositions
\ref{plane}, \ref{quadric}, \ref{cubic} we briefly describe all  possible
planes and  surfaces
of degree $2$ and $3$ contained in a rational normal $3$-fold $X\subset\PP^5$.
\end{remark}

\begin{enumerate}
\item
When $X=S(1, 1, 1)$. 
The only planes contained in $X$ are the ones in the ruling of $X$,
otherwise the linear system
$|\O_X(R)|$ would cut on a plane 
$\pi\subset X$, which is not in the ruling, a pencil of lines which of
course intersect eachother, while
the planes in $|\O_X(R)|$ are pairwise disjoint. 
The surfaces of degree $2$ contained in $X$ are either reducible and hence
linearly equivalent to $2R$, or irreducible and hence degenerate, therefore
linearly equivalent to $H-R$.
A reduced surface $Y\sim 2R$ is the disjoint union of two planes in $\PP^5$.
A surface $Q\sim H-R$ is a smooth quadric surface; 
the two systems of lines on $Q$
are cut by the linear system $|\O_X(R)|$ (lines of type $(1, 0)$ on $Q$)
 and by the linear system
$|\O_X(H-R)|$ (lines of type $(0, 1)$ on $Q$).
The surfaces of degree $3$ contained in $X$ are either reducible in the union of three 
planes and hence linearly equivalent to $3R$, or reducible in the union of 
an irreducible quadric surface and of a plane, hence  linearly equivalent
to $H-R+R\sim H$, or finally irreducible, therefore degenerate, and so linearly 
equivalent to an hyperplane section $H$.
A reduced surface $\sim 3R$ is the disjoint union of 
three planes in $\PP^5$.
A reducible hyperplane section of $X$ is the union of a smooth quadric surface 
and of a plane meeting along a line of  type $(1, 0)$.
An irreducible hyperplane section of $X$ is a smooth rational normal surface in $\PP^4$.
Lastly we recall that a surface $\sim aH+bR$ on $X$ is irreducible
 when $a=0$ and $b=1$, or $a>0$ and $b\geq -a$
(by \cite{h2}, V, 2.18, passing to general hyperplane sections).

\item
When $X=S(0, 1, 2)$, a plane contained in $X$ is either  one of the ruling
of $X$, therefore linearly equivalent to $R$, or  it is the plane $p\sim
H-2R$, i.e. the plane spanned by the vertex $V$ of $X$ and by the line 
image of the section defined by
$\operatorname{Proj}\O_{\PP^1}(1)\hookrightarrow \PP(\E)$. The reducible
surfaces of degree $2$ contained in $X$ are either linearly equivalent to
$2R$ (when reduced they are  the union of two planes meeting at the point
$V$),  or linearly equivalent to $H-R$
(the union of $p$ and of a plane $\pi\sim R$ meeting along a line passing through $V$),
 or linearly equivalent to
$2(H-2R)$ (the plane $p$ counted with multiplicity $2$). The  irreducible ones are
linearly equivalent to $H-R$.
An irreducible surface $Q\sim H-R$ is a quadric cone with vertex $V$; the pencil
 of lines on $Q$
is cut by the linear system $|\O_X(R)|$ (or equivalently by the linear system
$|\O_X(H-R)|$).
The surfaces of degree $3$ contained in $X$ are either reducible in the union of three 
planes and hence linearly equivalent to $3R$ (when reduced they are the union of
 three planes meeting at the point $V$),
 or to $2R+H-2R\sim H$ (when reduced each plane $\sim R$ meets the plane $p$ along a line passing through
the point $V$), or to 
$R+2(H-2R)=2H-3R$, or to
$3(H-2R)$ (in the last two cases the surface is not reduced). They
 may be also
 reducible in the union of 
an irreducible quadric cone and of a plane, hence linearly equivalent
to $H-R+R\sim H$ (the cone and the plane meet along a line passing through
$V$),
 or to $H-R+H-2R=2H-3R$ (the cone and the plane $P$ meet at the point
$V$).  Finally they can be irreducible, 
therefore degenerate and so linearly 
equivalent to an hyperplane section $H$, which is a rational normal surface in $\PP^4$.
As in the previous case, since a general hyperplane section of $X$ is smooth,
we have that a surface $\sim aH+bR$ on $X$ is irreducible
 when $a=0$ and $b=1$, or $a>0$ and $b\geq -a$.

\item
When $X=S(0, 0, 3)$ the situation is simpler since $\Cl(X)=\ZZ[R]$.
A plane contained in $X$ is a plane of the ruling. A surface of degree $2$ contained
in $X$ is always linearly equivalent to $2R$ and it is always reducible in the
union of two planes meeting along the line $l$ (the vertex of $X$). 
A surface of degree $3$ contained in $X$ is linearly equivalent
to $3R$ and it is reducible in the union of three planes meeting at $l$ 
if the proper transform is linearly
equivalent to $3\tilde R$, or it is irreducible, and hence a 
singular rational normal surface in $\PP^4$, if the proper transform is linearly
equivalent to $\tilde H$. 
In this case, since the a hyperplane section of $X$ is singular,
by \cite{h2}, V, 2.18 we have that a surface $S\subset X$ is irreducible 
if its proper transform
is $\tilde S\sim a \tilde H+b\tilde R$,
 with $a=0$ and $b=1$, or $a>0$ and $b\geq 0$.

\end{enumerate}

\begin{remark}
\label{description2}
Finally we want to describe how plane curves contained in $X$ look like.
Since $X$ is intersection of quadrics in $\PP^5$, then a plane curve 
of degree $\geq 3$ contained in $X$ must
lie in a plane of $X$. Therefore  we are interested just in lines and conics.
\end{remark}
\begin{enumerate}
\item
Let $X=S(1, 1, 1)$.
A line $r\subset X$ which is not contained in a plane $\pi\sim R$ of the scroll  is 
the base locus of a pencil of quadric surfaces $\sim H-R$, i.e. it is a line of type $(0, 1)$
(take the pencil of hyperplane sections passing through $r$ and a fixed plane 
$\sim R$ intersecting $r$).
A conic $C\subset X$ which does not lie on $\pi\sim R$ is contained in a quadric surface 
$Q\sim H-R$ (take a hyperplane section passing through $C$ and a plane
$\sim R$ meeting $C$) Therefore it is an hyperplane section of $Q$, i.e. a
curve of type $(1, 1)$ on $Q$.
\item
Let $X=S(0, 1, 2)$.
Every line $r\subset X$ is  contained in a plane 
 of the scroll. In fact if $r$ passes through $V$, then it is obviously
contained in some plane $\pi\sim R$. If $r$ does not pass through $V$,
then the
 plane spanned by $r$ and $V$ is contained in $X$. 
A conic $C\subset X$ 
which is not contained in a plane of the scroll and that
 passes through $V$ is reducible in the union of two lines.
If $C$ does not pass through $V$, then the cone over $C$ with vertex $V$
is a quadric cone $Q\sim H-R$, therefore $C$ is an hyperplane section of $Q$.
\item
Let $X=S(0, 0, 3)$.
A line $r\subset X$ is always contained in a plane 
 of the scroll $\pi\sim R$. In fact it is contained in the hyperplane section passing through
$r$ and the singular line $l$ of $X$, which splits in the union of three planes $\sim R$. 
A conic $C\subset X$ which does not lie on $\pi\sim R$ is 
 an hyperplane section of a surface $\sim 2R$, i.e. it is the union of two lines meeting
at a point.
\end{enumerate}
\section{Preliminaries}

We start by recalling a few results of \cite{ccd}.

From now on, let $C$ be an integral, nondegenerate curve of degree $d$ and
arithmetic genus $p_a(C)$ in $\PP^5$, with 
$d>\frac{2s}{3}\Pi_{i=1}^3 {(4!)}^{\frac{1}{4-i}}$. Assume $C$ not
contained on surfaces of degree $<s$ ($s\geq 4$) and define $m, \epsilon,
w, v, k,
\delta$ as in the statement of  Theorem \ref{teorema}.

Then the genus $p_a(C)$ is bounded by the function:
\[
G(d, 5, s)=1+\frac{d}{2}(m+w-2)-\frac{m+1}{2}(w-3)+\frac{vm}{2}(w+1)+\rho
\]
where $\rho=\frac{-\delta}{2}(w-\delta)$ if $\epsilon<w(4-v)$ and
$\rho=\frac{\epsilon}{2}-\frac{w}{2}(3-v)-\frac{\delta}{2}(w-\delta+1)$
if $\epsilon\ge w(4-v)$
(\cite{ccd}, section 5).

If $Z$ is a general hyperplane section of $C$ and $h_Z$ is the Hilbert
function of $Z$, then the difference $\Delta h_Z$ must be bigger than the
function $\Delta h$ defined by:
\[
\Delta h(n)=
\begin{cases}
0 & \text{if} \quad n<0 \cr
3n+1 & \text{if} \quad  0\le n \le w \cr
s & \text{if} \quad  w< n \le m \cr
s+k-3(n-m) & \text{if} \quad m< n \le m+\delta \cr
s+k-3(n-m)-1 & \text{if} \quad m+\delta < n \le m+w+e \cr
0 & \text{if} \quad n>m+w+e 
\end{cases}
\]
where $e=0$ if $\epsilon<w(4-v)$ and $e=1$ otherwise.

\begin{proposition}
\label{extcurve}
If $p_a(C)=G(d, 5, s)$, then $\Delta h_Z(n)=\Delta h(n)$ for all $n$ and
$C$ is arithmetically Cohen-Macaulay. Moreover $Z$ is contained on a
reduced curve $\Gamma$ of degree $s$ and maximal genus $G(s,
4)=\frac{w(w-1)(w-2)}{2}+wv$ in $\PP^4$ (Castelnuovo's curve). Since
$d>s^2$, $\Gamma$ is unique and, when we move the hyperplane, all these
curves $\Gamma$'s patch togheter giving a surface $S\subset \PP^5$ of
degree
$s$ through $C$.
\end{proposition}

\begin{pf} 
See \cite{ccd} 0.1, 6.1, 6.2, 6.3.
\end{pf}

S is a "Castelnuovo surface" in $\PP^5$, i.e. a surface whose general
hyperplane section is a curve of maximal genus in $\PP^4$.

\begin{proposition} 
S is irreducible and when $s\geq 9$ lies on a cubic rational
normal $3$-fold $X$ in
$\PP^5$ where it is cut by a hypersurface $G$ of degree $w+1$.
As a divisor on $X$ the surface $S$ is linearly equivalent to
$(w+1)H-(2-v)R$ (or $wH+R$ if $v=0$). 
\end{proposition}
\begin{pf}
$S$ is irreducible because $C$ is irreducible and is not contained on
surfaces of degree $<s$. 

A general hyperplane section $\Gamma$ of $S$ is a special Castelnuovo's
curve in
$\PP^4$ of degree $s$, then it lies on a rational normal cubic surface $W$
in $\PP^4$
which
is intersection of the quadric hypersurfaces containing $\Gamma$, hence
also $Z$; since $C$ is arithmetically Cohen-Macaulay these quadrics must
lift to quadric hypersurfaces in $\PP^5$ containing $C$, hence also $S$.
The intersection of these quadric hypersurfaces is a rational normal cubic
$3$-fold $X$
in $\PP^5$ whose general hyperplane section is $W$.

Moreover $\Gamma$ lies on a hypersurface of degree $w+1$ which does not
contain $W$; such a hypersurface must lift to a hypersurface $G$
of degree $w+1$ in $\PP^5$, containing $C$, hence containing $S$ since
$d>s^2$, and not containing $X$.
%
%
%
%
 \end{pf}

\begin{proposition}
There exists a hypersurface
$F$ of degree $m+1$, passing through $C$ and not containing $S$.
\end{proposition}
\begin{pf}
For a general hyperplane section $\Gamma$ of $S$, the Hilbert function
$h_\Gamma$ is known (see e.g. \cite{ha2}); in particular we have 
$\Delta
h_\Gamma (n)=\Delta h_Z (n)$ when $0\le n\le m$ and hence
$h^0({\I}_C (n))=h^0({\I}_S(n))$ when $0\le n\le m$. For $n=m+1$
one computes $\Delta h_Z (m+1)<\Delta h_\Gamma (m+1)$ and this implies
$h^0({\I}_C (m+1))>h^0({\I}_S(m+1))$.
\end{pf}

We recall here the definition of geometrical linkage.
\begin{definition}
Let $Y_1$, $Y_2$, $Y$ be subschemes of a projective space $\PP$, then
$Y_1$ and $Y_2$ are geometrically linked by $Y$ if
\begin{enumerate}
\item
$Y_1$ and $Y_2$ are equidimensional, have no embedded components and 
have no common components
\item $Y_1\cup Y_2=Y$, scheme theoretically.
\end{enumerate}
\end{definition}

\begin{definition}
\label{def1}
Call $C^\prime$ the curve residual to $C$ on  $S$ by $F$; $\deg C^\prime
=s-\epsilon-1$. Call $C^{\prime\prime}$ the curve residual to $C$ on $X$ by
$F$ and $G$.
\end{definition}

We note that $C^\prime$ is well defined since $S$ is irreducible and $F$
does
not contain $S$, $\deg C^\prime=s(m+1)-d=s-\epsilon-1$.
Moreover since $\deg C^\prime < \deg C$ the curve $C^\prime$ does not 
contain  $C$, which  is irreducible; therefore $C$ and $C^\prime$ {\it are
geometrically} linked by $S\cap F$.
 Also
$C^{\prime\prime}$ is well defined and $C^\prime\subset C^{\prime\prime}$:
\begin{itemize}
\item
if $s=3w+3$ ($v=2$), then
\[
C^{\prime\prime}=C^\prime.
\]
\item
if $s=3w+2$ ($v=1$), we can choose the plane $p_1\sim R$ linked to $S$ on
$X$ by
$G$ such that it is not contained in $F$, then
\[
C^{\prime\prime}=C^\prime+C_1
\]
where $C_1\subset p_1$ is a plane curve of degree $m+1$.
\item
if $s=3w+1$ ($v=0$) and $S\sim(w+1)H-2R$, we can choose the divisor $\sim
2R$
linked to $S$ on $X$ by $G$ such that it is the union of two distinct
plane $p_1$ and $p_2$ not contained in $F$, then
\[
C^{\prime\prime}=C^\prime+C_1+C_2
\]
where $C_1\subset p_1$ and $C_2\subset p_2$ are two planes curves of
degree $m+1$.
\item
if $s=3w+1$  and $S\sim wH+R$, we can choose the divisor $q\sim H-R$
linked to $S$ by $X$ and $G$ such that 
it is not contained in $F$, then
\[
C^{\prime\prime}=C^\prime+C_q
\]
where $C_q$  
is the intersection of $q$ and $F$.
\end{itemize}
$C$ and $C^{\prime\prime}$ are geometrically linked by $X\cap F\cap
G$
since they are equidimensional have no 
common components ($C$ is irreducible and $C^\prime$ does not contain $C$)
and with no embedded components
($X\cap  F\cap G$ is arithmetically Cohen Macaulay).

\begin{definitions}
\label{def2}
Let $X$ be a rational normal scroll in $\PP^5$. Call $A=Q_1\cap Q_2$ a generic complete 
intersection of two quadrics in $\PP^5$ containing $X$. Call $B\cong \PP^3$ the linked
scheme to $X$ by $A$. Call $Y_B=X\cap B$ the intersection scheme between $X$ and $B$.
Call $\ciii$ the curve linked to $C$ on $A$ by $F$ and $G$.
\end{definitions}

\begin{remark}
Hartshorne's Connectedness Theorem (\cite{e2}, Th. 18.12) implies that $Y_B$
is a divisor; moreover, since $Y_B\subset B$, it follows that $Y_B$
is contained in a hyperplane section and has degree $2$ (it is cut on $B$
by a quadric), therefore $Y_B\sim H-R$ as a divisor in $X$.
\end{remark}

\begin{definition}
\label{def3}
Call $Z$, $Z^\prime$, $Z^{\prime\prime}$, $\ziii$, $W$,  general hyperplane
sections
of
$C$, $C^\prime$, $C^{\prime\prime}$, $\ciii$, $X$ respectively.
By abusing notation call $A\subset\PP^4$ and $Y_B\subset W$ general hyperplane 
sections of $A\subset\PP^5$ and $Y_B\subset X$ respectively.
\end{definition}

Using the linkage techniques developed in \cite{f3} we 
can prove the following results
which are the main tool in the classification procedure.

\begin{lemma}
\label{first}
If $X$ is smooth or $X=S(0, 1, 2)$, then for $i\leq w$:
\begin{displaymath}
\begin{array}{rl}
h^0(\I_{C^{\prime\prime}|X}(iH+R)) \geq &
h^0(\I_{{\ziii}|W}(iH+R))=\\
= & h^1(\I_{Z|W}(m+w-i))=
\sum_{r=m+w-i+1}^\infty \Delta h(r).
\end{array}
\end{displaymath}
Moreover if $h^0(\I_{Z^{\prime\prime}|W}((i-1)H+R))=0$ and
$h^0(\I_{Z^{\prime\prime}|W}(iH+R))=h>0$, then 
$h^0(\I_{C^{\prime\prime}|X}((i-1)H+R))=0$
and $h^0(\I_{C^{\prime\prime}|X}(iH+R))=h$.
\end{lemma}
\begin{pf}
Since $C$ and $\cii$ are geometrically linked,
 $W$ is smooth and $C$
is arithmetically Cohen Macaulay, we apply \cite{f3}
Prop. 3.1 and we obtain:
$
h^0(\I_{Z^{\prime\prime}|W}(iH+R))= h^1(\I_{Z|W}(m+w-i))
=d-h_Z(m+w-i)
$. Then  note that for every $k$ 
we have $d-h_Z(k)=d+\Delta h_Z(k+1)- h_Z(k+1)=d+\sum_{r=k+1}^t \Delta
h_Z(r)-h_Z(t)
=\sum_{r=k+1}^\infty \Delta h_Z(r)$ because for $t$ big we
have $h_Z(t)=d$, and that, by Prop. \ref{extcurve},
$\Delta h_Z(r)=h (r)$.
By \cite{f3} Cor. 3.10 (if $X$ is smooth) or
Cor. 4.28 (if $X=S(0, 1, 2)$) we have that 
\begin{equation}
\label{bigger}
h^0(\I_{C^{\prime\prime}|X}(iH+R))\geq
h^0(\I_{Z^{\prime\prime}|W}(iH+R)). 
\end{equation}
From the exact sequence
\[
0\to \I_{C^{\prime\prime}|X}((k-1)H+R)\to
\I_{C^{\prime\prime}|X}(kH+R)\to
\I_{Z^{\prime\prime}|W}(kH+R)\to 0
\]
we obtain that if for $k=i-1$ we have
$h^0(\I_{Z^{\prime\prime}|W}(kH+R))=0$, then
$h^0(\I_{C^{\prime\prime}|X}(kH+R))=0$. In this
hypothesis for $k=i$ we have an injection
$H^0(\I_{C^{\prime\prime}|X}(iH+R)) \hookrightarrow 
H^0(\I_{Z^{\prime\prime}|W}(iH+R))$ and therefore 
by 
(\ref{bigger})
 $h^0(\I_{C^{\prime\prime}|X}(iH+R)) =
h^0(\I_{Z^{\prime\prime}|W}(iH+R))$.
\end{pf}

Using classical linkage techniques is easy to prove (\cite{f2} Lemma 4.7)
the following Lemma:
\begin{lemma}
\label{lemmac3}
Let $X=S(0, 0, 3)$, let $A$ as in Def. \ref{def2} and let $\ziii$ as in Def. \ref{def3}.
 Then for $i<w$
\[
h^0(\I_{{\ziii}|A}(i+1))=h^1(\I_{Z|X}(m+w-i)).
\]
\end{lemma}
Using Lemma \ref{lemmac3} and \cite{f3} Th. 4.20 we  prove the following
\begin{lemma}
\label{second}
If $X=S(0, 0, 3)$, then for $i\leq w$:
\[
h^0(\I_{C^{\prime\prime}|X}(iH+R))\geq
h^1(\I_{Z|W}(m+w-i))=\sum_{r=m+w-i+1}^\infty h(r).
\]
Moreover if $i<w$, 
 $h^0(\I_{\ziii|A}(i))=0$ and
$h^0(\I_{\ziii|A}(i+1))=h>0$, then 
$h^0(\I_{C^{\prime\prime}|X}((i-1)H+R))\geq h^0(\I_{\ciii|A}(i))= 0$ and
$h^0(\I_{C^{\prime\prime}|X}(iH+R))\geq h^0(\I_{\ciii|A}(i+1))=h$.
\end{lemma}
The next result is a formula which relates the arithmetic genera of the
curves $C$, $\cii$ and $Y=X\cap F\cap G$.

\begin{lemma}
Let $X$ be smooth or $X=S(0, 1, 2)$. Let $C$, $\cii$
and $Y=X\cap F\cap G$ as usual. Then we have the following relation:
\begin{equation}
\label{genus}
p_a(\cii)=p_a(C)-p_a(Y)+(m+w-1)\cdot\deg\cii+\deg(R\cap \cii)+1
\end{equation}
\end{lemma}
\begin{pf}
\cite{f3} Prop.
3.11 in the case $X=S(1, 1, 1)$ and  \cite{f3} Th. 4.30 if $X=S(0, 1, 2)$. 
\end{pf}

\section{The classification}

At this point we are able to prove the main Theorem.
The techniques that we use to prove Theorem \ref{teorema} are the same for the
three cases: 
$X=S(1, 1, 1)$,
$X=S(0, 1, 2)$ and
$X=S(0, 0, 3)$; therefore we  don't want to give a proof for all cases.
On the other side the analysis is slightly different case by case,
therefore, to be impartial, we will give the proof of Th. \ref{teorema}
2) in case $X=S(0, 1, 2)$, of Th. \ref{teorema} 3)
 in case $X=S(1, 1, 1)$
and of Th. \ref{teorema} 4)  in case $X=S(0, 0, 3)$. For a complete proof 
the reader may
consult \cite{f1}.
We will give a more precise description of such curves $\ci$ 
in the next propositions. The reader may go back to Remark
\ref{description1} where we have described planes and surfaces of degree
$2$ or
$3$  contained in
$X$, and to the previous section where we have introduced some notation.

\begin{pf}[of the Theorem \ref{teorema}]
\begin{enumerate}
\item
 Let $k=3$.
This happens if and only if $\epsilon=s-1$. It follows $\deg
C^\prime=s-\epsilon-1=0$ and we are done.
\item
 Let $k=2$ and let $X=S(0, 1, 2)$. We treat separately the cases 
$v=0, 1, 2$.

$\bullet$ Let $v=2$, i.e. $S\sim (w+1)H$. Then   $e=1$ and  we have
$\epsilon+1=2(w+1)+\delta$ with $0\leq \delta < w+1$.
By Lemma \ref{first}  we compute 
\[
h^0(\I_{{\ci}|X}(R))=1.
\]
Hence $\ci$ is contained in a plane $\pi\sim R$
and has degree $w+1\geq \deg\ci =3w+2-\epsilon\geq 1$.

$\bullet$ Let $v=1$, i.e. $S+p_1\sim (w+1)H$. If  $e=0$ (i.e. if 
$\epsilon < 3w$) we have
$\epsilon=2w+\delta$ with $0\leq 
\delta < w$. By Lemma \ref{first}  we compute 
$h^0(\I_{{\cii}|X}(R))=0$ and
$
h^0(\I_{{\cii}|X}(H+R))=3.
$
Since $(H+R)\cdot p_1\cdot H=1$
and $\deg C_1=m+1>1$, then all the surfaces $\sim H+R$ containing $\cii$
split 
in
the plane $p_1\supset C_1$ and in surfaces $\sim H-R$ containing $\ci$. 
Therefore
\[
h^0(\I_{{\ci}|X}(H))=3,
\] i.e. $\ci$ is contained in a plane $\pi$ and has degree
$w+1\geq\deg\ci=3w+1-\epsilon>1$. 
When $v=1$ we have  $e=1$ only if $\epsilon=3w$. In this case
$\deg\ci=1$, i.e. 
$\ci$ is a line. 

$\bullet$ Let $v=0$, then $e=0$ and we have $\epsilon=2w+\delta$
with
$0\leq 
\delta < w$.  By Lemma \ref{first} we compute 
$h^0(\I_{{\cii}|X}(R))=0$ and
$
h^0(\I_{{\cii}|X}(H+R))=2.
$
In case $S+p_1+p_2\sim (w+1)H$, 
since 
$(H+R)\cdot p_1\cdot H=(H+R)\cdot p_2\cdot H=1$
and $\deg C_1=\deg C_2=m+1>1$,
we find that
\[
h^0(\I_{{\ci}|X}(H-R))=2.
\]
If $\delta < w-1$, then $\deg\ci=3w-\epsilon=w-\delta>1$; since
$(H-R)\cdot (H-R)\cdot H=1$, then the linear system
$|\I_{{\ci}|X}(H-R)|$ has a fixed part which is necessarily the plane
$p\sim H-2R$.
If $\delta =w-1$, then $\ci$ is a line.
Therefore $\ci$ is contained in the plane $p\sim H-2R$ or it is a line.

In case $S+q\sim (w+1)H$, 
since 
$(H+R)\cdot q\cdot H=3$
and $\deg C_q=m+1>3$,
we find that
\[
h^0(\I_{{\ci}|X}(2R))=2.
\]
Therefore $\ci$ is contained in a plane $\pi\sim R$, which is the fixed part
of $|\I_{\ci|X}(2R)|$.

\item
 Let $k=1$ and let $X$ be smooth.

$\bullet$ Let $v=2$ and $e=0$; we have $\epsilon=w+\delta$ with $0\leq
\delta <w$.
By Lemma \ref{first}  we compute 
$h^0(\I_{{\ci}|X}(R))=0$
and
\[
h^0(\I_{{\ci}|X}(H+R))=3.
\]
Since  $\deg\ci=s-\epsilon-1\geq 5$ and $(H+R)\cdot (H+R)\cdot H=5$ we
deduce that the linear system $|\I_{{\ci}|X}(H+R)|$ has a fixed part which
has degree less or equal than $2$ as one can easily verify (if we suppose,
for example,  
that the fixed part is $L\sim H$, then
$h^0(\O_X(H+R-L))=h^0(\O_X(R))=2$, and we have a contradiction since
$h^0(\I_{{\ci}|X}(H+R-L))=3$). Therefore the fixed part can be of the
following types:
\begin{enumerate} \item
$\pi\sim R$. In this case $\ci$ is the union of a plane curve $\ci_1$ on
$\pi$ and
of a curve $\ci_2$ contained in the base locus of a net of hyperplane
sections,
i.e. in a plane $\sigma$. If $\sigma\sim R$, then the fixed part 
of $|\I_{{\ci}|X}(H+R)|$ is $\sim 2R$ and we are in the next case (b). The
other possibility is that $\sigma$ does not belong to $X$.
Since
$\deg\ci\geq w+3$ and $\pi\cdot S\cdot H=w+1$ this is possible only when
$\deg\ci=w+3$, i.e. $\delta=w-1$  and 
$\ci_2$ is a curve of type $(1, 1)$ on a quadric surface $\sim H-R$.
In this case $\ci$ is contained in the surface of degree two
$\pi\cup \sigma$.
\item
$Y\sim 2R$. Then 
$|\I_{{\ci}|X}(H+R)|=Y+|\O_X(H-R)|$ (in fact $h^0(\O_X(H-R))=3$). Since
$|\O_X(H-R)|$
is free, $\ci$ is contained in the surface of degree two $Y\sim 2R$.
\item
$Q\sim H-R$, 
i.e. $|\I_{{\ci}|X}(H+R)|=Q+|\O_X(2R)|$. Since $|\O_X(2R)|$ is free $\ci$
is contained in the smooth quadric surface $Q$.
\end{enumerate}
When $v=2$, then  $e=1$ only if $\epsilon=2w$, i.e. $\deg\ci=w+2$.
In this
case we write
$\epsilon+1=w+1+\delta$ with $\delta=w$.
By Lemma \ref{first} and  we compute 
$h^0(\I_{{\ci}|X}(R))=0$
and
\[
h^0(\I_{{\ci}|X}(H+R))=4.
\]
Since  $\deg\ci=w+2\geq 4$ and $(H+R)\cdot (H+R)\cdot H=5$ we deduce
that the linear system $|\I_{{\ci}|X}(H+R)|$ has a fixed part which
is, as one can easily verify, 
$\pi\sim R$. In this case $\ci$ is the union of a plane curve of degree
$w+1$ on
$\pi$ and of a line.
 
$\bullet$ Let $v=1$. Then $e=0$ and we have $\epsilon=w+\delta$ with
$0\leq\delta <w$.
By Lemma \ref{first} we compute 
$h^0(\I_{{\cii}|X}(R))=0$ and
$
h^0(\I_{{\cii}|X}(H+R))=2.
$
With the same computation we  have done previously (case $k=2$ and
$v=1$) one can deduce
that
\[
h^0(\I_{{\ci}|X}(H))=2.
\] 
Since $\deg\ci>w+1\geq 3$ and $H^3=3$ the linear pencil 
$|\I_{{\ci}|X}(H)|$ should have a fixed part, which can be of the
following types:
\begin{enumerate} 
\item
$\pi\sim R$. In this case $\ci$ is the union of a plane curve of degree
$w+1=\pi\cdot S\cdot H$
on
$\pi$ and
of 
a line, which is the base locus of a pencil of quadric surfaces $\sim
H-R$.
This is possible only when $\deg\ci=w+2$, i.e. $\delta=w-1$.
\item
$Q\sim H-R$, 
i.e. $|\I_{{\ci}|X}(H)|=Q+|\O_X(R)|$. Since $|\O_X(R)|$ is free $\ci$
is contained in the smooth quadric surface $Q$.
\end{enumerate}
The fixed part of $|\I_{\ci|X}(H)|$ cannot be $Y\sim 2R$ since in this
case we would have
$|\I_{{\ci}|X}(H)|=Y+|\O_X(H-2R)|$, while
$h^0(\O_{X}(H-2R))=0$.

$\bullet$ Let $v=0$, then  $e=0$ and we have $\epsilon=w+\delta$
with
$0\leq 
\delta < w$.  By Lemma \ref{first}  we compute 
$h^0(\I_{{\cii}|X}(R))=0$ and
$
h^0(\I_{{\cii}|X}(H+R))=1.
$
In case $S+p_1+p_2\sim (w+1)$, 
one easily deduces that
\[
h^0(\I_{{\ci}|X}(H-R))=1,
\] i.e. $\ci$ is contained in a smooth quadric
surface $Q\sim H-R$.
In case $S+q\sim (w+1)H$, 
one find that
\[
h^0(\I_{{\ci}|X}(2R))=1,
\]
therefore $\ci$ is contained in a reducible surface of degree two $Y\sim
2R$.

\item
 Let $k=0$ and let $X=S(0, 0, 3)$.
Then $e=0$ and we write $\epsilon=\delta$
with $0\leq \delta < w$.

$\bullet$ Let $v=2$.
By Lemma \ref{second}  we compute 
$h^0(\I_{{\ci}|X}(R))\geq h^0(\I_{{\ciii}|A}(1))=0$
and
\[
h^0(\I_{{\ci}|X}(4R))\geq h^0(\I_{{\ciii}|A}(2))=2.
\]
Since  $\deg\ci\geq 2w+3\geq 7$ and $4R\cdot 4R\cdot H=6$ by
(\ref{degree}), we
deduce that the linear system $|\I_{{\ci}|X}(4R)|$ has a fixed part.

We exclude that the fixed part is
$\pi\sim R$. Indeed in this case $\ci$ would be the union of a curve
contained in $\pi$ of degree at most  $\pi\cdot S\cdot H=w+1$,
and of a curve contained in the base locus of a pencil of 
hyperplane sections of degree at most 
$H^3=3$. But this is not possible  since $\deg\ci\geq w+5$.

Also
$Y\sim 2R$ is not possible. In fact in this case $\ci$
would be contained in  $Y$ (since the reduced singular line of
$X$, which is the base locus of 
the residual system $|\I_{\ci|X}(4R)-Y|\subset |\O_{X}(2R)|$, 
is contained in $Y$).
This can not happen: the quadric hypersurfaces which are  union of a
hyperplane
containing $Y$ and of a hyperplane containing $B$ 
 cut on $A$ a linear
subsystem of
$|\I_{\ciii|A}(2)|$ of projective dimension $2$ and this is not possible
since $h^0(\I_{{\ciii}|A}(2))=2$.

The only possibility is that the fixed part is $L\sim 3R$. In this 
case we claim that $\ci$ is contained in $L$.
In this case we prove that $\ci$ is contained in a surface of degree 
$3$ which is a 
hyperplane section of $X$.
$\ci$ is the union of a curve contained in $L$ and possibly of the reduced
singular line $l$ of $X$, which is the base locus of the residual system
$|\I_{\ci|X}(4R)-Y|=|\O_{X}(R)|$. If $L$ is reducible, then $l\subset L$
and we are done. If $L$ is irreducible, i.e. $L$ does not contain $l$,
then $\ci$ would contain $l$ with multiplicity $1$, but this is 
not possible, as the following computation shows.
Let 
$\tilde S\sim (w+1-a)\tilde H+3a\tilde R$ ($0\leq a \leq w$) and
$\tilde F\sim (m+1-b)\tilde H+3b\tilde R$ ($0\leq b\leq m$) 
be the proper transforms
 of $S$ and $F$ in $\tilde X$.
When $a\geq 1$ and $b\geq 1$
 by (\ref{intmult}) $\ci$ contains $l$ with multiplicity 
$
m(F, S; l)=3ab.
$

$\bullet$ Let $v=1$. 
By Lemma \ref{second} we compute 
$h^0(\I_{{\cii}|X}(R))\geq h^0(\I_{{\ciii}|A}(1))=0$
and
\[
h^0(\I_{{\cii}|X}(4R))\geq h^0(\I_{{\ciii}|A}(2))=1.
\]
Since by (\ref{degree}) $4R\cdot p_1\cdot H=2$ and $\deg C_1=m+1>2$, 
a surface 
$\sim 4R$ which contains $\cii$ splits in the union of $p_1\sim R$ and a
surface
$\sim 3R\sim H$ which contains $\ci$.
Therefore we have:
\[
h^0(\I_{{\ci}|X}(H))\geq 1.
\]
We claim that we have exactly
$
h^0(\I_{{\ci}|X}(H))= 1.
$ 
This follows from the fact that a quadric hypersurface that
is union of the hyperplane
$H$ which cuts $p_1+Y_B$ on $X$ and of a hyperplane which  cuts 
$L\in |\I_{{\ci}|X}(H)|$ is a quadric hypersurface which cuts on 
$A$ a divisor in 
$|\I_{{\ciii}|A}(2)|$ and we know that $h^0(\I_{{\ciii}|A}(2))=1$.

$\bullet$ Let $v=0$. By Lemma \ref{second}  we
compute 
$h^0(\I_{{\cii}|W}(4R))\geq h^0(\I_{{\ciii}|A}(2))=0$
and
\[
\begin{cases}
h^0(\I_{{\cii}|W}(7R))\geq h^0(\I_{{\ciii}|A}(3))=4
\qquad\text{if $\epsilon=w-1$}\\
h^0(\I_{{\cii}|W}(7R))\geq h^0(\I_{{\ciii}|A}(3))=3
\qquad\text{if $\epsilon<w-1$}\\
\end{cases}
\]
Since by (\ref{degree}) $7R\cdot p_1\cdot H=7R\cdot p_2\cdot H=3$ and $\deg C_1=\deg
C_2=m+1>3$ we deduce that
\[
\begin{cases}
h^0(\I_{{\ci}|W}(5R))\geq 4
\qquad\text{if $\epsilon=w-1$}\\
h^0(\I_{{\ci}|W}(5R))\geq 3
\qquad\text{if $\epsilon<w-1$}\\
\end{cases}
\]
and using similar arguments as
 in the previous case one can easily prove that equality holds.
We claim  that the linear system
$|\I_{{\ci}|W}(5R)|$ has a fixed part. To prove the claim we need first to analyze when $\ci$
may contain the singular line $l$ of $X$ as a component. Let 
$\tilde S\sim a\tilde H+(3w-3a+1)\tilde R$ ($0<a\leq w$) and
$\tilde F\sim (m+1-b)\tilde H+3b\tilde R$ ($0\leq b\leq m$) be the proper transforms
 of $S$ and $F$ in $\tilde X$.
When $b\geq 1$ by (\ref{intmult}) $\ci$ contains $l$ with multiplicity 
$
m(F, S; l)=3b(w-a)+b\geq b.
$
On the other hand $\ci$ is contained in the scheme 
$S\cap D$ for some $D\in |\I_{{\ci}|W}(5R)|$; 
since $S\cap D$ contains $l$ with multiplicity
$m(D, S; l)=2(w-a)+1$ if $\tilde D\sim \tilde H+2\tilde R$  or 
$m(D, S; l)=5(w-a)+2$ if $\tilde D\sim 5\tilde R$, then $m(F, S; l)$
should be less or equal then these values.
Therefore when $b\geq 1$, since $3b(w-a)+b>2(w-a)+1$, we exclude the
possibility
 $\tilde D\sim \tilde H +2\tilde R$; when $b\geq 2$, since $3b(w-a)+b>5(w-a)+2$, we exclude
the possibility $\tilde D\sim 5 \tilde R$.
Hence $\ci$ may contain $l$ only if $b=1$, i.e. $\tilde F\sim m\tilde H+3\tilde R$,
and the divisors in the linear system
$|\I_{\ci|X}(5R)|$ are all reducible in the union of five planes.
In this case $\ci$ contains $l$ with multiplicity $m(F, S; l)=3(w-a)+1$,
which has to be less or equal than $m(F, 5R; l)$, that is $5$ by
(\ref{intmult}). This is possible only if either $a=w-1$, or $a=w$.
When $a=w-1$ we have $\tilde S\sim (w-1)\tilde H+4\tilde R$ and $\ci$ contain $l$
with multiplicity $4$.
When $a=w$ we have $\tilde S\sim w\tilde H+\tilde R$ and $\ci$ contain $l$
with multiplicity $1$.

Now we are able to prove that  $|\I_{\ci|X}(5R)|$  has a fixed part.
Let us suppose first that $\ci$ contains $l$. 
In this case if 
 the linear system  $|\I_{\ci|X}(5R)|$ has no fixed part, then $\ci$ 
is supported on 
the line $l$. Our previous computation
implies that $\deg\ci\leq 4$, while we know that
$\deg\ci\geq 7$. Therefore $|\I_{\ci|X}(5R)|$ has a fixed part as claimed.
Let us suppose that $\ci$ does not contain $l$. If 
$|\I_{\ci|X}(5R)|$ has no fixed part,
then the generic element $D$ in $|\I_{\ci|X}(5R)|$ is irreducible and has
proper transform
$\tilde D\sim \tilde H+2\tilde R$, therefore for $D, D^\prime$ in the
linear system we have
$m(D, D^\prime; l)=1$. In this case, since by (\ref{degree}) 
$5R\cdot 5R\cdot H=8$, 
the base locus of a pencil in 
$|\I_{\ci|X}(5R)|$ not supported on $l$ has degree $7$. Since
$\deg\ci\geq 7$ 
and $h^0(\I_{\ci|X}(5R))=3$  we have a contradiction.
Therefore $|\I_{\ci|X}(5R)|$ has a fixed part.

We claim first that this fixed part cannot be
$\pi\sim R$. In this case $\ci$ would be the union of a curve
$\ci_1\subset\pi$
and of a curve $\ci_2$ contained in the base locus of a linear subsystem 
$|\alpha|\subset |\I_{\ci|X}(4R)|$ of projective dimension $3$ 
if $\epsilon=w-1$, $2$
if $\epsilon<w-1$. We want to prove that $|\alpha|$ has a fixed part.
Let us suppose first that $\ci$ contains $l$, with multiplicity $4$ (if $a=w-1$) or $1$
(if $a=w$)
 by our previous
computation. 
If $|\alpha|$ has no fixed part then $\ci_2$ is supported on
$l$, and
since the component of $\ci_1$ disjoint from the line $l$ has degree 
equal to $\tilde\pi\cdot\tilde S\cdot\tilde H= a$, we should have 
$\deg\ci=4+w-1=w+3$ (if $a=w-1$) or $\deg\ci=1+w$ (if $a=w$), but this is not possible 
since $\deg\ci=s-\epsilon-1=3w-\epsilon>2w\geq w+3$.
 Therefore we have a contradiction and $|\alpha|$ has a fixed part.
With similar arguments it is easy to prove that $|\alpha|$ has a fixed part
if we suppose that $\ci$ does not contain $l$.
The fixed part of $|\I_{\ci|X}(5R)|$ can be of the following types:
\begin{enumerate}
\item
$Y\sim 2R$. In this case $\ci=\ci_1\cup\ci_2$, where $\ci_1\subset Y$ and
$\ci_2$ is contained in the base locus of a linear system
$|\beta|\subset |\O_X(H)|$ of projective dimension $3$ if $\epsilon=w-1$,
$2$ if $\epsilon<w-1$.
If $\epsilon=w-1$, then $\ci_2$ is a line $r\subset\pi\sim R$. Therefore
$\ci$ is contained in the cubic surface $Y\cup\pi\sim 3R$.
If $\epsilon<w-1$, then $\ci_2\subset \sigma$ 
is a plane curve contained in
a plane $\sigma$. 
If $\ci$ contains $l$, the divisors in $|\beta|$ are
reducible in the union of $3$ planes, therefore $\sigma\sim R$ is a fixed
part for
$|\I_{\ci|X}(5R)|$ and we are in the next case (b). 
If $\ci$ does not contain $l$, then
$\deg \ci_1\leq 2w$; therefore if $\epsilon\leq w-3$ we have
$\deg\ci_2\geq 
3$ which implies $\sigma\sim R$. If $\epsilon=w-2$ and $\deg\ci_1=2w$ 
(i.e. $\tilde S\sim w\tilde H+\tilde R$), then
$\deg\ci_2=2$ and $\sigma$ may be a plane not contained in $X$. $\ci$ is
contained in the cubic surface $Y\cup \sigma$.
\item
$L\sim 3R$. Here we must put $\epsilon<w-1$, since for $\epsilon=w-1$ we
have
$h^0(\I_{\ci|X}(5R))=4$, while
$h^0(\O_{X}(5R-L))=
h^0(\O_{X}(2R))=3$. In this case $\ci$ is contained in a cubic surface
$L\sim 3R$, hyperplane section of $X$.
 \end{enumerate}
\end{enumerate}
\end{pf}

In the next propositions we give a closer description 
of case 2) in Th. \ref{teorema} 
when $X=S(0, 1, 2)$,
of case 3) when $X$ is smooth and of case 4) when $X=S(0, 0, 3)$. In
all the other cases 
one can give a  similar description, but we don't intend
to go  through this.

\begin{proposition}
\label{plane}
Let $k=2$. If $X=S(0, 1, 2)$ then we have the
following possibilities:
\begin{enumerate}
\item
When $v=2$ then $\ci$ is a plane curve of degree $1\leq\deg\ci\leq w+1$
contained in a plane $\pi\sim R$. $\ci$ does not pass through the singular point 
$V$ of
$X$.
\item
When $v=1$ then $\ci$ is  
contained in a plane which can be of the following types:
\begin{enumerate}
\item
if $\deg\ci=w+1$, i.e. $\epsilon=2w$, then $\ci\subset \pi\sim R$ and
passes through the vertex $V$ of $X$.
\item
if $\deg\ci\leq w$, i.e. $2w < \epsilon \leq 3w$, then $\ci$ lies either in
$\pi\sim R$
or in $p\sim H-2R$ and may pass or not through $V$.
\item
if $\deg\ci=2$, i.e. if $\epsilon=3w-1$ there is the further 
possibility that $\ci$ is a conic lying on a 
plane $\sigma$ which does not belong to the scroll.
In this case $\ci$ is either the union of two lines passing through $V$,
or it is a hyperplane section of a 
quadric cone $\sim H-R$.
\end{enumerate}
\item
When $v=0$, $S\sim (w+1)H-2R$ 
and $1\leq \deg\ci\leq w-1$ (i.e. $2w<\epsilon\leq 3w-1$) then $\ci$ is contained in
$p\sim H-2R$.
$\ci$ may pass or not  through the
vertex of $X$.
When $\deg\ci=1$, i.e. if $\epsilon=3w-1$, there is the further
possibility that $\ci$ is a line contained in a plane $\pi\sim R$ and
passing through $V$. 
\item
When $v=0$ and $S\sim wH+R$, then $\ci$ is a plane curve
of degree $1\leq\deg\ci\leq w$ contained in $\pi\sim R$.
$\ci$ may pass or not through the vertex of $X$.
\end{enumerate}
\end{proposition}
\begin{pf}
\begin{enumerate}
\item
This is the case $k=v=2$ in the proof of the Theorem \ref{teorema}.
That $\ci$ does not pass through the point  $V$
follows  by genus formula (\ref{genus}): in this
case in fact $\deg(R\cap\ci)\geq 1$ and this would  imply that $p_a(C)$
is not
maximal.
\item
Looking at the case $k=2$ and $v=1$ in the proof of  Theorem \ref{teorema}
we know
that $\ci$ is contained in a plane $\pi$ which is the base locus of a net
of
hyperplanes. If $\deg\ci\geq 3$ this plane must be contained in the
scroll, therefore is 
either $\pi\sim R$ or $\pi=p\sim
H-2R$.
\begin{enumerate}
\item
When $\deg\ci=w+1$, since $S\cdot p \cdot H=w$ we must have $\pi\sim R$.
Since $S\sim (w+1)H-R$, the hypersurface $G$ of degree $w+1$ which cut $S$
on
$X$ must pass through the point $V$; therefore the curve $\ci$ must
contain $V$.  
\item
When $\deg\ci < w+1$, by genus formula (\ref{genus}) one verifies  
that $\pi$ can be both $\pi\sim R$ or $\pi=p\sim H-2R$ and that $\ci$
may pass or not through the vertex $V$.
\item
If $\deg\ci= 2$ the plane that contains $\ci$ may not belong to the
scroll. Therefore, as we have explained 
in Remark \ref{description2},
either $\ci$ is the union of
two lines meeting at $V$,  or it is a hyperplane section of a quadric cone
$\sim H-R$
\end{enumerate}
\item
Looking at the proof of  Theorem \ref{teorema} when $k=2$ and $S\sim (w+1)H-2R$
we know that if $\deg\ci>1$ then $\ci$ lies on the plane $p\sim H-2R$.
In this case, since $S\cdot p\cdot H=w-1$, 
we must exclude the  case
$\deg\ci=w$. 
If $\deg\ci=1$ then $\ci$ is a line which can also be the base locus of a
pencil of quadric cones $\sim H-R$, therefore in this case it lies in a
plane $\pi\sim R$ and passes through $V$.
\item
This is proved in  Theorem \ref{teorema}. 
Since $\deg\ci=3w-\epsilon$ and $2w\leq
\epsilon\leq 3w-1$, we have that $1\leq\deg\ci\leq w$.
\end{enumerate}
\end{pf}

\begin{proposition}
\label{quadric}
If $k=1$, then $\ci$ is contained in a surface of degree two.
When $X=S(1, 1, 1)$
we have the
following possibilities:
\begin{enumerate}
\item
if $S \nsim wH+R$ and $\epsilon=w, w+1$ the surface may be a smooth
quadric
surface $Q\sim H-R$. In this case if $\epsilon=w$ then $\ci$ is a curve of
kind $(w-1+v, w+1)$ on $Q$; if $\epsilon=w+1$ then $\ci$ is of kind $(
w-1+v, w)$ on $Q$.
\item
If $v=2$ or $v=0$ and $S\sim wH+R$ the surface may be reducible in the
union of two disjoint planes $\pi_1\sim R$ and $\pi_2\sim R$. In this case
$\ci=\ci_1\cup\ci_2$, where $\ci_1$ is a curve of degree $2w+1-\epsilon$
if $v=2$ (resp. of degree $2w-\epsilon$ if $S\sim wH+R$) on $\pi_1$ and
$\ci_2$ is a curve of degree $w+1$ (resp. $w$) on $\pi_2$.
\item{} 
When a) $v=2$ and $\epsilon=2w-1$ or b) $v=2$ and $\epsilon=w$ or c)
$v=1$
and
$\epsilon=2w-1$, the surface may be the union of a plane $\pi\sim R$
and a plane $\sigma$ not contained in the scroll $X$.
In this cases $\ci=\ci_1\cup\ci_2$, where $\ci\subset\pi$ has degree $w+1$
and $\ci_2\subset \sigma$ is a) a conic of kind $(1, 1)$ or b), c)
a line of kind $(0, 1)$. $\ci_1$ and $\ci_2$ intersect each other in a
point.
\item
if $v=1$ and $\epsilon\neq w, w+1, 2w-1$ or  $S\sim (w+1)H-2R$
and
$\epsilon\neq w, w+1$ there are no curves of maximal genus $G(5, d, s)$ on
$X=S(0, 0, 1)$.
 \end{enumerate}
\end{proposition}
\begin{pf}
\begin{enumerate}
\item
Let $S$ and $\Gamma$ be as usual.
Looking at the cases $v=2$ part c), $v=1$ part b)
and $S\sim (w+1)H-2R$ in the proof of Theorem \ref{teorema}
we know that
$\ci$ may be contained in a smooth quadric surface $Q\sim  H-R$. Hence
the hypersurfaces in $\PP^4$ of degree $m+1$ through $Z$ cut on $\Gamma$
a linear series $|\gamma|$ of degree $\deg\zi=s-\epsilon-1$, whose
divisors
are contained in  conics $\sim H-R$, and therefore in  hyperplane sections.
Since $\Gamma$ is arithmetically Cohen Macaulay, it follows that $|\gamma|$
is cut out on $\Gamma$ by the linear system of hyperplanes through 
$H\cdot\Gamma-\deg\zi=\epsilon+1$ points.
The projective dimension of $|\gamma|$ is equal to
$h_\Gamma(m+1)-h_Z(m+1)$. If $\epsilon\geq w+1$ we have
$\dim |\gamma|=1$, therefore the $\epsilon +1$ points should span a
$\PP^2$.
If $\epsilon=w$ we have $\dim |\gamma|=2$, hence they lie on a line
$l\sim R$ on $W$.
If $\epsilon\geq w+1$ the movable part of $|\gamma|$ is cut out by the
pencil of
conics
$\sim H-R$ through a fixed point (since $h^0(\O_W(H-R))=3$); this implies
that  the remaining $\epsilon$ points lie on a line $\sim R$. Therefore
$\epsilon+1=\Gamma\cdot R+1=w+2$, i.e. $\epsilon=w+1$.
In this case $\zi$ is a set of $\Gamma\cdot (H-R)-1=2w+v-1$ points on the
conic $\sim H-R$, hence $\ci$ is linked to a line of kind $(0, 1)$ on a
smooth quadric surface $Q\sim H-R$
by the intersection with $S$, i.e. $\ci$ is a curve of type $(w-1+v, w)$.
If $\epsilon=w$, the movable part of $|\gamma|$ is cut out by the linear
system 
$|\O_W(H-R)|$ and $\zi$ is a set of $\Gamma\cdot (H-R)=2w+v$ points on a
conic $\sim H-R$. Therefore $\ci$ is the intersection of $S$ with a smooth
quadric surface $Q$, hence is a curve of type $(w-1+v, w+1)$.
\item
We omit the proof in this case, since it goes 
on similarly as the previous one. 
\item
This third possibility follows directly from the
proof of  Theorem \ref{teorema}.
\item
If $v=1$ from the proof of  Theorem \ref{teorema} we know that $\ci$ lies on a smooth
quadric surface $\sim H-R$ or, if $\epsilon=2w-1$, in the surface $\pi\cup
\sigma$ of part 3) in the statement.
Therefore by part 1) of this proposition we deduce that for $\epsilon\neq
w, w+1, 2w-1$ the curve $\ci$ does not exist.
In a similar way one deduces that $\ci$ does not exist when $S\sim
(w+1)H-2R$ and $\epsilon\neq w, w+1$.
\end{enumerate}
\end{pf}

The proof of the next proposition  is omitted since the argument is
similar to the previous one. Also the explicit computation for the
multiplicity of the line $l$ as a component of the linked curve $\ci$ is
similar to the one appearing in the proof of the Theorem part 4).

\begin{proposition}
\label{cubic}
If $k=0$, then $\ci$ is contained in a surface of degree $3$. When
$X=S(0, 0, 3)$ we have the following possibilities:
\begin{enumerate}
\item
The surface is a hyperplane section $L$ of $X$. More precisely:
\begin{enumerate}
\item
if $\epsilon=0$, then $\ci$ is linked to a line by the intersection
$S\cap L$; 
\item
if $\epsilon=1$, then $\ci$ is linked to a conic by the intersection
$S\cap L$; 
\item
if $\epsilon>1$, then $L$ split in the union of three planes $\pi_i\sim R$
($i=1, 2, 3$). In this case, when $v=2$ the proper transform of $S$
is $\tilde S\sim (w+1-a)\tilde H +3a\tilde R$ with $0\leq a \leq
w-\epsilon$. $\ci$ is the union of the line $l$ counted with multiplicity
$3a$, of a curve $\ci_1\subset \pi_1$ of degree $w+1-a$, of a curve
$\ci_2\subset\pi_2$ of the same degree, and of a curve
$\ci_3\subset \pi_3$ of degree $w-\epsilon-a$.
$\ci_1$ and $\ci_2$ intersect each other in $w+1-a$ points on $l$. $\ci_1$
and
$\ci_2$ both intersect $\ci_3$ in $w-\epsilon-a$ points on $l$.

When $v=1$ the proper transform of $S$
is $\tilde S\sim a\tilde H +(3w-3a+2)\tilde R$ with $\epsilon+1\leq a \leq
w$. $\ci$ is the union of the line $l$ counted with multiplicity
$3(w-a)+2$, of a curve $\ci_1\subset \pi_1$ of degree $a$, of a curve
$\ci_2\subset\pi_2$ of the same degree, and of a curve
$\ci_3\subset \pi_3$ of degree $a-\epsilon-1$.
$\ci_1$ and $\ci_2$ intersect each other in $w+1-a$ points on $l$. $\ci_1$
and
$\ci_2$ both intersect $\ci_3$ in $w-\epsilon-a$ points on $l$.

When $v=0$ the proper transform of $S$
is $\tilde S\sim w\tilde H +\tilde R$.
If $\epsilon=w-1$, then $\ci$ does not contain the line $l$ and it is the
union of a curve 
$\ci_1\subset
\pi_1$ of degree $w$, of a curve
$\ci_2\subset\pi_2$ of the same degree, and of a curve
$\ci_3\subset \pi_3$ of degree $w-\epsilon$.
$\ci_1$ and $\ci_2$ intersect each other in $w$ points on $l$. $\ci_1$
and
$\ci_2$ both intersect $\ci_3$ in $w-\epsilon$ points on $l$.
If $\epsilon<w-1$, then $\ci$ may contain the line $l$
 with multiplicity
$1$. In this case $\ci$ is the union of $l$, 
$\ci_1$, of 
$\ci_2$, and of a curve
$\ci_3\subset \pi_3$ of degree $w-\epsilon-1$.
$\ci_1$ and $\ci_2$ intersect each other in $w$ points on $l$. $\ci_1$
and
$\ci_2$ both intersect $\ci_3$ in $w-\epsilon-1$ points on $l$.
\end{enumerate}
\item
If $v=0$ and $\epsilon=w-2$ the surface may split in the union of a plane
$\pi_1\sim R$, of a plane $\pi_2\sim R$ and of a plane $\sigma$ which is
not contained in the scroll. 
$\ci$ does not contain the line $l$ as a component and  is the union of
a curve
$\ci_1\subset
\pi_1$ of degree $w$, of a curve
$\ci_2\subset\pi_2$ of the same degree, and of the lines $r_1$ and $r_2$
of intersection between $\sigma$ and $X$.
$\ci_1$ and $\ci_2$ intersect eachother in $w$ points on $l$. $r_1$
and
$r_2$  intersect $\ci_1$ and $\ci_2$ in theyr common point
on $l$.
 \end{enumerate}
\end{proposition}

\section{The existence}

Lastly we need an effective construction for curves of degree $d$, 
genus $G(d, 5, s)$ in $\PP^5$, not lying on surfaces of degree $<s$.
It should be noted that in case $k=v=1$ it is not possible to construct
curves of maximal genus on a {\it smooth} rational normal $3$-fold 
(see Prop. \ref{quadric} case 4)). In this case the construction is possible
only on a rational normal $3$-fold whose vertex is a point.
Before the construction we state the following, easy to prove, 
result (see \cite{ro} Lemma 1 pg. 133), that we will use.

\begin{lemma}
\label{rogora}
Let $X$ be a smooth $3$-fold. Let $\Sigma$ be a linear system of surfaces
of $X$ and let $\gamma$ be a curve contained in the base locus of $\Sigma$.
Suppose that the generic surface of $\Sigma$ is smooth at the generic 
point of $\gamma$ and that it has at least a singular point which is variable 
in $\gamma$. Then all the surfaces of $\Sigma$ are tangent along $\gamma$.
\end{lemma}
\begin{example}
For all $d>\frac{2s}{3}\Pi_{i=1}^3 {(4!)}^{\frac{1}{4-i}}$ and $s\geq 4$, 
there is a smooth curve $C\subset \PP^5$ of degree $d$, genus $G(d, 5, s)$
contained on an irreducible surface $S$ of degree $s$.
\end{example}

$\bullet$ Let $k=3$. In this case take $S$ to be the complete intersection
of a smooth rational normal $3$-fold $X\subset \PP^5$ 
and a general hypersurface $G$ of degree $w+1$
in $\PP^5$. The complete intersection of $S$ with a general hypersurface $F$ of
degree $m+1$  gives the required curve $C$.
\medskip

$\bullet$ Let $k=2$ and $v=2$. In this case we have
$\epsilon+1=2(w+1)+\delta$ with $0\leq \delta\leq w$.
Let $\pi$ be a plane contained in a smooth rational normal $3$-fold
$X$ in $\PP^5$ and
let $D\subset\pi$ be a smooth plane curve of degree \,
$w+1>\deg D=\epsilon-2w-1\geq 0$ (possibly $D=\emptyset$). 
Let us consider the linear system $|\I_{D|X}(w+1)|$, which is not empty since
 contains 
the linear subsystem $L+|\O_X(w)|$, where $L$ is a hyperplane section of
$X$ containing the plane $\pi$. This shows also that 
$|\I_{D|X}(w+1)|$ is not composed with a pencil 
because in this case every element in the system
would be a sum of algebraically equivalent divisors, but, for example, 
the elements in $L+|\O_X(w)|$ are obviously 
not of this type. 
Since $\deg D<w+1$ the linear system $|\I_{D|\pi}(w+1)|$ is not
empty and
its base locus is the curve $D$; this shows that the base locus of
$|\I_{D|X}(w+1)|$ is exactly the curve $D$, because $|\O_X(w)|$
is base points free.
Therefore by Bertini's Theorem the general divisor in 
$|\I_{D|X}(w+1)|$ is an irreducible surface $S$ of degree $s$, 
which is smooth outside  $D$. 
We claim that  $S$ is in fact smooth at every point $p$ of $D$.
To see this, by Lemma \ref{rogora}, it is enough to prove that,
for every point $p\in D$, there exists 
a surface in
 $|\I_{D|X}(w+1)|$ which is smooth at $p$, and that 
for a generic point  $q\in D$, there exist two surfaces in 
$|\I_{D|X}(w+1)|$  with distinct 
tangent planes at $q$.
Indeed, for every $p\in D$ we can always find a  surface $T$ in the 
linear system $|\O_X(w)|$ which does not pass through $p$, therefore
the surface $L+T$ is smooth at $p$ with tangent plane $\pi$. Moreover
a generic surface in the linear system $|\I_{D|X}(w+1)|$ which cut $D$
on $\pi$ has at $p$ tangent plane $T_p\neq \pi$.

Let $\ci\subset\pi$ be the linked curve to $D$ by the intersection of $S$
and $\pi$.
Let us consider the linear system $|\I_{\ci|S}(m+1)|$. Since 
$\deg\ci=s-\epsilon-1<m+1$, with the same arguments
used above it is easy to see that $|\I_{\ci|S}(m+1)|$
is not empty,
is not composed with a  pencil
and has base
locus equal to  the curve $\ci$.
Therefore by Bertini's Theorem
 the generic curve $C$ in the movable part of
this linear system is irreducible and smooth. Moreover $C$ lies on a
smooth surface $S$ of degree $s=3w+3$ and it has the required
numerical characters, as one may easily verify 
using in the genus
formula (\ref{genus})
$
p_a(\ci)=\frac{1}{2}(3w+2-\epsilon-1)(3w+2-\epsilon-2)
$
(computed by Clebsch's formula) and $\deg (R\cap \ci)=0$.
\medskip 

$\bullet$ Let $k=2$ and $v=1$. In this case we have
$\epsilon=2w+\delta$ with $0\leq \delta\leq w$. Let $X$ and $\pi\subset X$ be
as in the previous case and let $p_1\neq \pi$ be 
an other plane contained in $X$.
Let $D\subset\pi$ be a smooth plane curve of degree 
\, $0\leq\deg D=\epsilon-2w<w+1$
on $\pi$.
 In this case by Bertini's Theorem we can find an irreducible
 surface 
$S\sim (w+1)H-R$ of degree $s=3w+2$ in the movable part of the 
linear system $|\I_{D\cup p_1|X}(w+1)|$, smooth outside $D$.
 With
the same argument used in the previous case we prove that $S\cup p_1$ is 
smooth at every
point of $D$. Namely, for every $p\in D$,
a generic surface in the linear system $L+L_1+|\O_X(w-1)|$, where
$L=\pi+Q$ and
$L_1=p_1+Q_1$ are reducible hyperplane sections containing respectively
$\pi$ and $p_1$ and
such that $p\notin Q\cup Q_1$, is smooth at $p$ with tangent
plane 
$\pi$, while a surface in the linear system $L_1+|\I_{D|X}(w)|$ which cut
$D$ on $\pi$, has tangent plane $T_p\neq \pi$. 
Since $D\cap p_1=\emptyset$
this implies that $S$ is smooth.

Let $\ci$ be the linked
curve to $D$ by the intersection $\pi\cap S$. 
By Bertini's Theorem the linked curve $C$ to $\ci$ by the intersection of
$S$ with a general element
$F_{m+1}$ in the linear system $|\I_{\ci|S}(m+1)|$ is smooth of degree $d$.
Moreover it 
lies on a smooth surface $S$ of degree $s=3w+2$ and it has 
the required genus, as one 
can compute by using formula (\ref{genus}), where the curve $\cii$
is the union of $\ci$ and of the curve $C_1\subset p_1$ cut on $p_1$ by 
$F_{m+1}$,
and $\deg (R\cap\cii)=0$.

\medskip
$\bullet$ Let $k=2$ and $v=0$. In this case we have
$\epsilon=2w+\delta$ with $0\leq \delta <w$. Let $X$ and $\pi\subset X$
as before and let $q$ be a smooth quadric surface contained in $X$, 
intersecting $\pi$ along a line $r$.
Let $D\subset\pi$ be a smooth plane curve of degree 
\, $0\leq\deg D=\epsilon-2w<w$
on $\pi$ (when $\deg D=1$ we suppose that $D$ does not coincide with the line 
$r$). In this case by Bertini's Theorem we can find an irreducible
 surface 
$S\sim wH+R$ of degree $s=3w+1$ in the movable part of the 
linear system $|\I_{D\cup q|X}(w+1)|$, smooth outside $D$. 
As in the previous cases we claim that $S$ is smooth.
Namely, for every $p\in D$, we can find in the  movable part of our linear
system
$|\I_{D\cup q}(w+1)|$,   
a surface 
 which is smooth at $p$ with tangent plane equal to $\pi$
(take a surface of the form
$L+\pi^\prime+T$, with $L$ as usual, $\pi^\prime\sim R$ a plane of $X$
disjoint from $\pi$
and $T$ a surface in $|\O_X(w-1)H|$ that does not pass through $p$), and 
a surface
with tangent plane $T_p\neq \pi$ ( take a surface
of the form $\pi^\prime+V$, where $V$ is a
divisor in $\I_{D|X}(w)$ which cut
$D$ on $\pi$).

Let $\ci$ be the linked
curve to $D$ by the intersection $\pi\cap S$; we have $\deg\ci=3w-\epsilon$. 
The linked curve $C$ to $\ci$ by the intersection of $S$ with a general element
$F_{m+1}$ in the linear system $|\I_{\ci|S}(m+1)|$ is smooth of degree $d$ and lies
on the smooth surface $S$ of degree $s=3w+1$. 
Its genus is maximal as one  can compute by
 using formula (\ref{genus}), where the curve $\cii$
is the union of $\ci$ and of the curve $C_q\subset q$ of type $(m+1, m+1)$ on $q$,
intersecting each other in $3w-\epsilon$ points along $r$. In this case
$\deg (R\cap\cii)=m+1$.

\medskip

$\bullet$ Let $k=1$ and $v=2$. In this case we have 
$\epsilon=w+\delta$ with 
$0\leq\delta\leq w$.  
Let $\pi_1$ be a plane contained in a smooth rational normal $3$-fold
$X$ in $\PP^5$ and
let $D\subset\pi_1$ be a smooth plane curve of degree \, $0\leq\deg
D=\epsilon-w<w+1$ (possibly $D=\emptyset$). 
In this case by Bertini's Theorem we can find an irreducible
 surface 
$S$ of degree $s=3w+3$ in the 
linear system $|\I_{D|X}(w+1)|$, which is also smooth (by using
exactly the same argument used in case $k=v=2$).
Let $\pi_2$ be an other plane contained in $X$ but not contained in $S$ 
and let $\ci$ be the linked curve to $D$ by the intersection $S\cap(\pi_1\cup\pi_2)$.
Therefore $\ci$ is the disjoint union of a plane curve of degree $2w+1-\epsilon$
on $\pi_1$ and of a plane curve of degree $w+1$ on $\pi_2$.
Let us consider the linear system $|\I_{\ci|S}(m+1)|$, whose base
locus is the curve $\ci$.
Therefore by Bertini's Theorem
 the generic curve $C$ in the movable part of
this linear system is  smooth, of degree $d$, lies on a smooth surface $S$
of degree
$s=3w+3$  and it has the required
genus as one may easily verify   using genus formula (\ref{genus}).

\medskip

$\bullet$ Let $k=1$ and $v=1$. In this case $\epsilon=w+\delta$ with $0\leq \delta < w$.
Let $X$ be a rational normal $3$-fold singular at a point $V$, i.e. $X=S(0, 1, 2)$.
Let $p_1\sim R$ be a plane in the
ruling of $X$ intersecting the plane $p\sim H-2R$ along a line $r_1$. 
Let $D\subset p$ be a smooth plane curve of degree 
\, $0\leq \deg D=\epsilon-w<w$ contained in the plane
$p$ and not passing through the vertex $V$ of $X$ (in particular when $\deg D=1$,
the curve $D$ can not be the line $r_1$).
By Bertini's Theorem we can find an irreducible 
surface $S\sim (w+1)H-R$  of degree
$s=3w+2$ in the movable part of the linear system $|\I_{D\cup p_1|X}(w+1)|$ 
smooth outside $D$ and $V$. 
With the same argument used in case $k=2$ and $v=1$ 
one can
say that
$S$ is smooth outside the vertex $V$
of $X$.
Let $\pi\sim R$ be an other plane in the ruling of $X$ and not 
contained in $S$; let $\ci$
be the linked curve to $D$ by the intersection $S\cap (p\cup\pi)$.  
Therefore $\ci$ is the  union of a plane curve $\ci_p$ of degree $2w-\epsilon$
on $p$ which does not pass
through $V$, and of a plane curve $\ci_\pi$ of degree $w+1$ on $\pi$ 
passing once through $V$, that
meet each other in $2w-\epsilon$ points
on the line $r$ of intersection between $p$ and $\pi$.
Let us consider the linear system $|\I_{\ci|S}(m+1)|$, whose base
locus is the curve $\ci$.
Therefore by Bertini's Theorem
 the generic curve $C$ in the movable part of
this linear system is irreducible, of degree $d$, lies on an irreducible
surface $S$ of degree $s=3w+2$ and it is smooth outside $V$.
But since $V$ does impose just one condition to the linear system  
$|\I_{\ci|S}(m+1)|$, 
namely the one imposed by $\ci_\pi$, it follows that $C$ does not pass through $V$,
hence $C$ is smooth.
The genus of $C$ is maximal and  can be computed by using formula
(\ref{genus}),  where the curve $\cii$
is the union of $\ci$ and of a plane curve $C_1\subset p_1$ of degree $m+1$
intersecting along $2w-\epsilon$ points on the line $r_1=\pi\cap p_1$,  
and $\deg (R\cap\cii)=2w-\epsilon+1$.
\medskip

$\bullet$ Let $k=1$ and $v=0$. In this case $\epsilon=w+\delta$ with $0\leq \delta < w$.
Let $X$ be again a smooth rational normal $3$-fold in $\PP^5$.
Let $q\sim H-R$ be a smooth quadric surface contained in $X$ and let $\pi_1\sim R$
be a plane of $X$ intersecting $q$ along a line $r_1$.
Let $D\subset \pi_1$ be a smooth plane curve of degree 
\, $0\leq \deg D=\epsilon-w<w$ contained in $\pi_1$ (when $\deg D=1$ we
suppose that $D$ does not coincide with the line $r_1$).
With the same argument used in case $k=2$ and $v=0$ 
we  find an irreducible smooth surface $S\sim wH+R$ of
degree
$s=3w+2$ in the movable part of the linear system $|\I_{D\cup q|X}(w+1)|$. 
 Let $\pi_2\sim R$ be an other plane  of $X$ not contained in
$S$; let $\ci$ be the linked curve to $D$ by the intersection $S\cap
(\pi_1\cup\pi_2)$.   Therefore $\ci$ is the disjoint union of a plane curve
of degree $2w-\epsilon$ on $\pi_1$ and of a plane curve of degree $w$ on
$\pi_2$. Let us consider the linear system $|\I_{\ci|S}(m+1)|$, whose base
locus is the curve $\ci$.
Therefore by Bertini's Theorem
 the generic curve $C$ in the movable part of
this linear system is smooth, of degree $d$ and lies on a smooth surface $S$ of degree $s=3w+1$.
The genus of $C$ is maximal and can be computed using formula (\ref{genus}), 
where the curve $\cii$
is the union of $\ci$ and of a curve $C_q\subset q$ of type $(m+1, m+1)$ on $q$,
intersecting at $2w-\epsilon$ points on the line $r_1=\pi_1\cap q$ and at $w$ points  
on the line $r_2=\pi_2\cap q$. Moreover  
$\deg (R\cap\cii)=m+1$.

\medskip
$\bullet$
Let $k=0$ and $v=2$. In this case we have that $0\leq \epsilon<w$.
Let $X$ be a smooth rational normal $3$-fold in $\PP^5$
and let $\pi$ be a plane of $X$. Let $D\subset\pi$ be a plane curve of degree
$0\leq\deg D=\epsilon+1<w+1$
in $\pi$.
With the same argument used in case $k=v=2$ one can prove that the general divisor in 
the linear system $|\I_{D|X}(w+1)|$ is a smooth surface $S$ of degree $s=3w+3$.
Let $L$ be a general hyperplane section of $X$ containing $\pi$,
i.e. the union of a smooth quadric surface $Q$ and $\pi$ meeting 
along a line $r$.
Let $\ci$ be the linked curve to $D$ by the intersection $S\cap L$, i.e. $\ci$
is the union of a curve $\ci_Q\subset Q$ of type $(w+1, w+1)$ and of a plane 
curve
$\ci_\pi\subset \pi$ of degree $w-\epsilon$, meeting at $w-\epsilon$ points.
Let us consider the linear system $|\I_{\ci|S}(m+1)|$. By Bertini's Theorem
 the generic curve $C$ in the movable part of
this linear system is smooth of degree $d$. Moreover $C$ lies on a
smooth surface $S$ of degree $s=3w+3$ and it has the required
numerical characters, as one may easily verify 
using
$
p_a(\ci)=w^2+{\frac{1}{2}}(w-\epsilon-1)(w-\epsilon-2)+w-\epsilon-1
$
(computed by Noether's formula) and $\deg (R\cap \ci)=w+1$ in the genus
formula (\ref{genus}).
\medskip

$\bullet$ Let $k=0$ and $v=1$. Again we have $0\leq\epsilon <w$.
Let $X$, $\pi\subset X$ and $p_1$ as in the case 
$k=2$ and $v=1$.
Let $D\subset \pi$ be a smooth
plane curve of degree $\epsilon+1$. 
With
the same argument used in the  case $k=2$ and $v=1$ 
we find an irreducible smooth surface
$S$ in the movable part of a general divisor in 
the linear system
$|\I_{D\cup p_1|X}(w+1)|$.
Let
$L=\pi\cup Q$ be a general hyperplane section of $X$ containing the plane $\pi$, 
which intersects
the plane $p_1$ in a line $r_1$ of type $(1, 0)$.
$\ci$ is the linked curve to $D$ by the intersection $S\cap L$;
therefore
$\ci$ is the union of a curve $\ci_Q\subset Q$ of type $(w, w+1)$ and of
a plane 
curve
$\ci_\pi\subset \pi$ of degree $w-\epsilon$, meeting at $w-\epsilon$ points.
The generic curve $C$ in the movable part of
this linear system 
$|\I_{\ci|S}(m+1)|$ 
is smooth of degree $d$. Moreover $C$ lies on a
smooth surface $S$ of degree $s=3w+2$ and it has the required
numerical characters as one may easily verify using 
formula (\ref{genus}). In this case $\cii$ is the union of $\ci$ with a plane
curve of degree $m+1$ on $p_1$ meeting at $w+1$ points, and $\deg (R\cap \cii)=w+1$.
\medskip

$\bullet$ Let $k=0$ and $v=0$. We have $0\leq\epsilon <w$.
 Let $X$, $\pi\subset X$ and
$L=\pi\cup Q$ as in the previous case. Let $q\sim H-R$ be an other 
smooth quadric surface contained in $X$ which intersects $Q$ along a line $r$ of type
$(0, 1)$ and $\pi$ along a line $s$ of type $(1, 0)$. Let $D\subset \pi$ be a smooth
plane curve of degree $\epsilon$ (when $\deg D=1$ we suppose that $D$ does not 
coincide with the line $s$). 
With the same argument used in  case  $k=2$ and $v=0$ we find an
irreducible  smooth
 surface $S$ is the movable part of a general divisor in the linear system
$|\I_{D\cup q|X}(w+1)|$.
$\ci$ is the linked curve to $D$ by the intersection $S\cap L$;
therefore
$\ci$ is the union of a curve $\ci_Q\subset Q$ of type $(w+1, w)$ and of
a plane 
curve
$\ci_\pi\subset \pi$ of degree $w-\epsilon-1$, meeting at $w-\epsilon-1$ points.
The generic curve $C$ in the movable part of
this linear system 
$|\I_{\ci|S}(m+1)|$ 
is smooth of degree $d$. Moreover $C$ lies on a
smooth surface $S$ of degree $s=3w+1$ and it has the required
numerical characters as one may easily verify by using 
formula (\ref{genus}). In this case $\cii$ is the union of $\ci$ with a
curve of type $(m+1, m+1)$ on $q$ meeting at $w+1$ points along the 
line $r$ and at $w-\epsilon-1$ points along the line $s$, and $\deg (R\cap \cii)=w+m+1$.


\begin{thebibliography}{XX}

\bibitem[CCD]{ccd} L. Chiantini, C. Ciliberto, V. Di Gennaro:
{\em The genus of projective curves}.
Duke Math. J., 50, 760, 2, 1993, 229-245.

\bibitem[CC]{cc} L. Chiantini, C. Ciliberto:
{\em Curves of maximal genus in $ \PP^4$}. Proceedings
"Zero-dimensional-schemes, Ravello 1992", De Gruyter, 1994. 

\bibitem[dC]{dc} M. A. A de Cataldo: {\em The genus of curves on the three
dimensional quadric}. Nagoya Math. J., 147, 1997, 193-211.

\bibitem[E1]{e1} D. Eisenbud: {\em Linear sections of determinantal
varieties}.
Amer. J. Math., 110, 1992, 541-575.

\bibitem[E2]{e2} D. Eisenbud: {\em Commutative Algebra with a View toward
Algebraic Geometry}. GTM 150, Springer Verlag, New York, 1994.

\bibitem[EH]{eh} D. Eisenbud, J. Harris: 
{\em On varieties of minimal degree ( a centennial account), Proceedings
of
the
AMS Summer Institute in Algebraic Geometry, Bowdoin, 1985}. Proceedings of
Symposia in pure Math., 46, AMS, 1987.

\bibitem[F1]{f1} R. Ferraro: {\em Curve di genere massimo in $ \PP^5$ and
Explicit Resolutions of Double Point Singularities of Surfaces}. Tesi di
Dottorato in Matematica, IX ciclo, Universit\`a di Roma "Tor Vergata",
1998. Unpublished work.

\bibitem[F2]{f2} R. Ferraro: {\em Weil divisors on rational normal
scrolls.}
Geometric and Combinatorial Aspects of Commutative Algebra.
Herzog-Restuccia Editors.
Lectures Notes in Applied
and Pure Mathematics, Marcel Dekker Ser. Vol. 217, pp. 
183-198, 2001. (http://arXiv.org/abs/math.AG/0105081). 

\bibitem[F3]{f3} R. Ferraro: {\em Linkage on singular rational normal
surfaces and thre-folds with application to the classification of curves
of maximal genus}. Preprint  (http://arXiv.org/abs/math.AG/0105094).


\bibitem[GP]{gp} L. Gruson, C. Peskine: {\em Genre des courbes dans
l'espace 
projectif}. Springer Lect. Notes 687, 1978.



\bibitem[Ha1]{ha1} J. Harris: {\em A bound on the geometric genus of
projective
varieties}. Ann. Scuola Norm. Pisa, 8, 1981, 35-68. 

\bibitem[Ha2]{ha2} J. Harris: {\em Curves in projective space}. Presses de
l'Universite de Montreal, 1982. 

\bibitem[H1]{h1} R. Hartshorne: {\em Generalized divisors on Gorenstein 
schemes.} K-Theory Journal, 8, 1994, pp. 287-339.

\bibitem[H2]{h2} R. Hartshorne: {\em Algebraic Geometry}. 
GTM 52, Springer-Verlag, New York, 1977. 

\bibitem[PS]{ps} C. Peskine, L. Szpiro:
{\em Liaison des vari\'et\'es alg\'ebriques (I)}. Invent. Math., 26, 1974,
271-302.

\bibitem[R]{r} A. P. Rao: {\em Liaison among curves in $\PP^3$}. Invent.
Math., 50, 1979, 205-217.

\bibitem[Ro]{ro} E. Rogora: {\em Metodi proiettivi e differenziali 
per lo studio di alcune questioni relative alle variet\`a immerse}.
Tesi di Dottorato, Universit\`a degli Studi di Roma ``La Sapienza'', 1996.
Unpublished work. 

\bibitem[S]{s} F. O. Schreyer: {\em Syzygies of canonical curves and
special
linear series}, Math. Ann., 275, 1979, 105-137.

\end{thebibliography}
\end{document}